\documentclass[11pt]{article}
\usepackage{cite}
\usepackage{mathrsfs}
\usepackage{amssymb,amsfonts,amsmath}
\usepackage{enumerate}
\usepackage{indentfirst}
\usepackage{latexsym,bm}
\usepackage[noend]{algpseudocode}
\usepackage{algorithmicx,algorithm}
\usepackage{algorithm}
\usepackage{makeidx}
\usepackage{fancybox}
\usepackage{color}
\usepackage{multicol}
\usepackage[latin1]{inputenc}
\usepackage{graphicx}
\usepackage{pstricks,pst-node,pst-text,pst-3d}
\usepackage{tikz}
\newtheorem{thm}{Theorem}[section]

\newtheorem{lem}[thm]{Lemma}
\newtheorem{op}[thm]{Problem}

\newtheorem{pro}[thm]{Proposition}

\newtheorem{obs}[thm]{Observation}
\newtheorem{conj}{Conjecture}[section]

\newenvironment{pf}{{\noindent \it \bf Proof:}}{{\hfill$\Box$}\\}

\begin{document}

\title{\bf Extremal results for directed tree connectivity}
\author{Yuefang Sun$^{1}$\thanks{Supported by Zhejiang Provincial
Natural Science Foundation (No. LY20A010013).}\qquad\\
$^{1}$E-mail: yuefangsun2013@163.com,\\
School of Mathematics and Statistics, Ningbo University,\\
Zhejiang 315211, P. R. China\\
}
\maketitle

\begin{abstract} For a digraph $D=(V(D), A(D))$, and a set $S\subseteq V(D)$ with
$r\in S$ and $|S|\geq 2$, an $(S, r)$-tree is an out-tree $T$ rooted at $r$ with
$S\subseteq V(T)$. Two $(S, r)$-trees
$T_1$ and $T_2$ are said to be arc-disjoint if $A(T_1)\cap
A(T_2)=\emptyset$. Two arc-disjoint $(S, r)$-trees $T_1$ and $T_2$ are
said to be internally disjoint if $V(T_1)\cap V(T_2)=S$. Let $\kappa_{S,r}(D)$ and $\lambda_{S,r}(D)$ be the maximum
number of internally disjoint and arc-disjoint $(S, r)$-trees in $D$,
respectively. The generalized $k$-vertex-strong connectivity
of $D$ is defined as
$$\kappa_k(D)= \min \{\kappa_{S,r}(D)\mid S\subset V(D), |S|=k, r\in
S\}.$$ Similarly, the generalized $k$-arc-strong connectivity
of $D$ is defined as
$$\lambda_k(D)= \min \{\lambda_{S,r}(D)\mid S\subset V(D), |S|=k, r\in
S\}.$$
The generalized $k$-vertex-strong connectivity and generalized $k$-arc-strong
connectivity are also called directed tree connectivity which could be seen as a generalization of classical connectivity of digraphs.

A digraph $D=(V(D), A(D))$ is called minimally generalized $(k,
\ell)$-vertex (respectively, arc)-strongly connected if $\kappa_k(D)\geq \ell$ (respectively, $\lambda_k(D)\geq \ell$) but for any arc $e\in A(D)$, $\kappa_k(D-e)\leq \ell-1$ (respectively, $\lambda_k(D-e)\leq \ell-1$). In this paper, we study the minimally generalized $(k, \ell)$-vertex (respectively, arc)-strongly connected digraphs. We compute the minimum and maximum sizes of these digraphs, and give characterizations of such digraphs for some pairs of $k$ and $\ell$. 
\vspace{0.3cm}
\vspace{0.3cm}

\end{abstract}


\section{Introduction}\label{sec:intro}

We refer the readers to \cite{Bang-Jensen-Gutin, Bang-Jensen-Gutin2,
Bondy} for graph theoretical notation and terminology not given
here. \footnote{
Note that all digraphs considered in this paper have no parallel
arcs or loops.} For a graph $G=(V,E)$ and a set $S\subseteq V$ of at
least two vertices, an {\em $S$-Steiner tree} or, simply, an {\em
$S$-tree} is a tree $T$ of $G$ such that
$S\subseteq V(T)$. Two $S$-trees $T_1$ and $T_2$ are said to be {\em
edge-disjoint} if $E(T_1)\cap E(T_2)=\emptyset$. Two edge-disjoint
$S$-trees $T_1$ and $T_2$ are said to be {\em internally disjoint}
if $V(T_1)\cap V(T_2)=S$. The basic problem of {\sc Steiner Tree
Packing} is defined as follows: the input consists of an undirected
graph $G$ and a subset of vertices $S\subseteq V(D)$, the goal is to
find a largest collection of edge-disjoint $S$-Steiner trees.
Besides this classical version, people also study some other
variations, such as packing internally disjoint Steiner trees,
packing directed Steiner trees and packing strong subgraphs \cite{Cheriyan-Salavatipour,
DeVos-McDonald-Pivotto, Kriesell, Lau, Sun-Gutin-Ai, Sun-Gutin-Yeo-Zhang, Sun-Yeo,  West-Wu}.

An {\em out-tree (respectively, in-tree)} is an oriented tree in which every vertex except one, called the {\em root}, has in-degree (respectively, out-degree) one. An {\em out-branching (respectively, in-branching)} of $D$ is a spanning out-tree (respectively, in-tree) in $D$.
For a digraph $D=(V(D), A(D))$, and a set $S\subseteq V(D)$ with
$r\in S$ and $|S|\geq 2$, a {\em directed $(S, r)$-Steiner tree} or,
simply, an {\em $(S, r)$-tree} is an out-tree $T$ rooted at $r$ with
$S\subseteq V(T)$ \cite{Cheriyan-Salavatipour}. Two $(S, r)$-trees
$T_1$ and $T_2$ are said to be {\em arc-disjoint} if $A(T_1)\cap
A(T_2)=\emptyset$. Two arc-disjoint $(S,r)$-trees $T_1$ and $T_2$ are
said to be {\em internally disjoint} if $V(T_1)\cap V(T_2)=S$.

Cheriyan and Salavatipour \cite{Cheriyan-Salavatipour} introduced
and studied the following two type of directed Steiner tree packing
problems. {\sc Arc-disjoint directed Steiner tree packing} (ADSTP):
The input consists of a digraph $D$ and a subset of vertices
$S\subseteq V(D)$ with a root $r$, the goal is to find a largest
collection of arc-disjoint $(S, r)$-trees. {\sc Internally-disjoint
directed Steiner tree packing} (IDSTP): The input consists of a
digraph $D$ and a subset of vertices $S\subseteq V(D)$ with a root
$r$, the goal is to find a largest collection of internally disjoint
$(S, r)$-trees.

Sun and Yeo \cite{Sun-Yeo} introduced the concept of directed tree connectivity which is
related to directed Steiner tree packing problems and extends the concept of tree connectivity of undirected graphs to digraphs (see e.g. \cite{Hager, Li-Mao5, Li-Mao-Sun}). Let $\kappa_{S,r}(D)$ and $\lambda_{S,r}(D)$ be the maximum
number of internally disjoint and arc-disjoint $(S, r)$-trees in $D$,
respectively. The {\em generalized $k$-vertex-strong connectivity}
of $D$ is defined as
$$\kappa_k(D)= \min \{\kappa_{S,r}(D)\mid S\subset V(D), |S|=k, r\in
S\}.$$ Similarly, the {\em generalized $k$-arc-strong connectivity}
of $D$ is defined as
$$\lambda_k(D)= \min \{\lambda_{S,r}(D)\mid S\subset V(D), |S|=k, r\in
S\}.$$ By definition, $\kappa_2(D)=\kappa(D)$ and $\lambda_2(D)=\lambda(D)$. Therefore, these two parameters could be seen as generalizations of vertex-strong connectivity and arc-strong connectivity of a digraph, respectively. The generalized $k$-vertex-strong connectivity and generalized $k$-arc-strong connectivity are also called {\em directed tree connectivity}.

Now we introduce more concepts and parameters related to directed tree connectivity. A digraph $D=(V(D), A(D))$ is called {\em minimally generalized $(k,
\ell)$-vertex (respectively, arc)-strongly connected} if $\kappa_k(D)\geq \ell$ (respectively, $\lambda_k(D)\geq \ell$) but for any arc $e\in A(D)$, $\kappa_k(D-e)\leq \ell-1$ (respectively, $\lambda_k(D-e)\leq \ell-1$).
By the definition of $\kappa_k(D)$ (respectively, $\lambda_k(D)$) and Theorem \ref{thma}, we clearly have $2\leq k\leq n, 1\leq \ell \leq n-1$.

Let $\mathfrak{F}(n,k,\ell)$ be the set of all minimally generalized $(k,
\ell)$-vertex-strongly connected digraphs with order $n$. We define
$$F(n,k,\ell)=\max\{|A(D)| \mid D\in \mathfrak{F}(n,k,\ell)\}$$ and
$$f(n,k,\ell)=\min\{|A(D)| \mid D\in \mathfrak{F}(n,k,\ell)\}.$$ We
further define $$Ex(n,k,\ell)=\{D\mid D\in \mathfrak{F}(n,k,\ell),
|A(D)|=F(n,k,\ell)\}$$ and $$ex(n,k,\ell)=\{D\mid D\in
\mathfrak{F}(n,k,\ell), |A(D)|=f(n,k,\ell)\}.$$

Similarly, let $\mathfrak{G}(n,k,\ell)$
be the set of all minimally generalized $(k,
\ell)$-arc-strongly connected
digraphs with order $n$. We define
$$G(n,k,\ell)=\max\{|A(D)| \mid D\in \mathfrak{G}(n,k,\ell)\}$$ and
$$g(n,k,\ell)=\min\{|A(D)| \mid D\in \mathfrak{G}(n,k,\ell)\}.$$
We further define $$Ex'(n,k,\ell)=\{D\mid D\in
\mathfrak{G}(n,k,\ell), |A(D)|=G(n,k,\ell)\}$$ and
$$ex'(n,k,\ell)=\{D\mid D\in \mathfrak{G}(n,k,\ell),
|A(D)|=g(n,k,\ell)\}.$$
By definition, we directly have $f(n,k,\ell)\leq F(n,k,\ell)$ and $g(n,k,\ell)\leq G(n,k,\ell)$.

The complexity for the decision versions of {\sc IDSTP} and {\sc ADSTP} on general digraphs, Eulerian digraphs and symmetric digraphs have been completely determined, as shown in the following tables.

\begin{center}
\begin{tabular}{|c||c|c|c|} \hline
\multicolumn{4}{|c|}{Table 1: Directed graphs} \\ \hline
$\kappa_{S,r}(D) \geq \ell$? & $k=3$                                    & $k \geq 4$     &  $k$ part \\
$ |S|=k$                     &                                          & constant       &  of input \\ \hline  \hline
$\ell =2$                    & NP-complete \cite{Cheriyan-Salavatipour} & NP-complete\cite{Sun-Yeo}    &  NP-complete\cite{Sun-Yeo}  \\ \hline
$\ell \geq 3$ constant       & NP-complete \cite{Sun-Yeo}                             & NP-complete\cite{Sun-Yeo}    &  NP-complete\cite{Sun-Yeo}  \\ \hline
$\ell$ part of input         & NP-complete \cite{Sun-Yeo}                             & NP-complete\cite{Sun-Yeo}    &  NP-complete \cite{Sun-Yeo} \\ \hline
\end{tabular}
\end{center}

\begin{center}
\begin{tabular}{|c||c|c|c|} \hline
\multicolumn{4}{|c|}{Table 2: Directed graphs} \\ \hline
$\lambda_{S,r}(D) \geq \ell$? & $k=3$                                    & $k \geq 4$     &  $k$ part \\
$ |S|=k$                     &                                          & constant       &  of input \\ \hline  \hline
$\ell =2$                    & NP-complete \cite{Cheriyan-Salavatipour} & NP-complete\cite{Sun-Yeo}    &  NP-complete\cite{Sun-Yeo}  \\ \hline
$\ell \geq 3$ constant       & NP-complete \cite{Sun-Yeo}                             & NP-complete \cite{Sun-Yeo}   &  NP-complete\cite{Sun-Yeo}  \\ \hline
$\ell$ part of input         & NP-complete \cite{Sun-Yeo}                             & NP-complete\cite{Sun-Yeo}    &  NP-complete\cite{Sun-Yeo}  \\ \hline
\end{tabular}
\end{center}

\begin{center}
\begin{tabular}{|c||c|c|c|} \hline
\multicolumn{4}{|c|}{Table 3: Eulerian digraphs} \\ \hline
$\kappa_{S,r}(D) \geq \ell$? & $k=3$        & $k \geq 4$     &  $k$ part \\
$ |S|=k$                      &              & constant       &  of input \\ \hline  \hline
$\ell =2$                     & NP-complete \cite{Sun-Yeo}            & NP-complete  \cite{Sun-Yeo}             &  NP-complete\cite{Sun-Yeo}   \\ \hline
$\ell \geq 3$ constant        & NP-complete  \cite{Sun-Yeo}           & NP-complete  \cite{Sun-Yeo}             &  NP-complete \cite{Sun-Yeo}  \\ \hline
$\ell$ part of input          & NP-complete \cite{Sun-Yeo} & NP-complete \cite{Sun-Yeo}   &  NP-complete \cite{Sun-Yeo}  \\ \hline
\end{tabular}
\end{center}

\begin{center}
\begin{tabular}{|c||c|c|c|} \hline
\multicolumn{4}{|c|}{Table 4: Eulerian digraphs} \\ \hline
$\lambda_{S,r}(D) \geq \ell$? & $k=3$        & $k \geq 4$     &  $k$ part \\
$ |S|=k$                      &              & constant       &  of input \\ \hline  \hline
$\ell =2$                     & Polynomial \cite{Sun-Yeo}  & Polynomial\cite{Sun-Yeo}     &  Polynomial\cite{Sun-Yeo}   \\ \hline
$\ell \geq 3$ constant        & Polynomial\cite{Sun-Yeo}   & Polynomial\cite{Sun-Yeo}     &  Polynomial\cite{Sun-Yeo}   \\ \hline
$\ell$ part of input          & Polynomial\cite{Sun-Yeo}   & Polynomial\cite{Sun-Yeo}     &  Polynomial\cite{Sun-Yeo}   \\ \hline
\end{tabular}
\end{center}

\begin{center}
\begin{tabular}{|c||c|c|c|} \hline
\multicolumn{4}{|c|}{Table 5: Symmetric digraphs} \\ \hline
$\kappa_{S,r}(D) \geq \ell$? & $k=3$        & $k \geq 4$     &  $k$ part \\
$ |S|=k$                     &              & constant       &  of input \\ \hline  \hline
$\ell =2$                    & Polynomial \cite{Sun-Yeo}  & Polynomial \cite{Sun-Yeo}    &  NP-complete \cite{Sun-Yeo} \\ \hline
$\ell \geq 3$ constant       & Polynomial \cite{Sun-Yeo}  & Polynomial \cite{Sun-Yeo}    &  NP-complete \cite{Sun-Yeo} \\ \hline
$\ell$ part of input         & NP-complete \cite{Sun-Yeo} & NP-complete \cite{Sun-Yeo}   &  NP-complete \cite{Sun-Yeo} \\ \hline
\end{tabular}
\end{center}

\begin{center}
\begin{tabular}{|c||c|c|c|} \hline
\multicolumn{4}{|c|}{Table 6: Symmetric digraphs} \\ \hline
$\lambda_{S,r}(D) \geq \ell$? & $k=3$        & $k \geq 4$     &  $k$ part \\
$ |S|=k$                      &              & constant       &  of input \\ \hline  \hline
$\ell =2$                     & Polynomial \cite{Sun-Yeo}  & Polynomial\cite{Sun-Yeo}     &  Polynomial\cite{Sun-Yeo}   \\ \hline
$\ell \geq 3$ constant        & Polynomial \cite{Sun-Yeo}  & Polynomial \cite{Sun-Yeo}    &  Polynomial\cite{Sun-Yeo}   \\ \hline
$\ell$ part of input          & Polynomial \cite{Sun-Yeo}  & Polynomial\cite{Sun-Yeo}     &  Polynomial\cite{Sun-Yeo}   \\ \hline
\end{tabular}
\end{center}


For directed tree connectivity, some inequalities concerning
parameters $\kappa_k(D)$ and $\lambda_k(D)$ were obtained. Let $D$
be a strong digraph with order $n$. For $2\leq k\leq n$, It was proved in \cite{Sun-Yeo}
that $1\leq \kappa_k(D)\leq n-1$ and $1\leq \lambda_k(D)\leq n-1$. Moreover, all bounds are sharp, and
those digraphs $D$ for which $\kappa_k(D)$~(respectively, $\lambda_k(D))$ attains
the upper bound are characterized. Then the authors studied the relation between the
directed tree connectivity and classical connectivity of digraphs by
showing that $\kappa_k(D)\leq \kappa(D)$ and $\lambda_k(D)\leq
\lambda(D)$. Furthermore, these bounds are sharp. In
the same paper, the sharp Nordhaus-Gaddum type bounds for
$\lambda_k(D)$ were also obtained; moreover, extremal digraphs for the lower bounds were characterized. 

In this paper, we will study the minimally generalized $(k,
\ell)$-vertex-strongly connected digraphs and minimally generalized $(k,
\ell)$-arc-strongly connected digraphs. We first give
characterizations of such digraphs for some pairs of $k$ and $\ell$ (Theorem~\ref{thmi}), and then obtain exact values or sharp bounds
for the functions $f(n,k,\ell)$, $F(n,k,\ell)$, $g(n,k,\ell)$ and $G(n,k,\ell)$ (Theorem~\ref{thmh}). Some open problems will also be posed.

\paragraph{Additional Terminology and Notation.} For a digraph $D$, its {\em reverse} $D^{\rm rev}$ is a digraph with the same vertex set such that
$xy\in A(D^{\rm rev})$ if and only if $yx\in A(D)$. A digraph $D$ is {\em symmetric} if $D^{\rm rev}=D$. In other words, a symmetric digraph $D$ can be obtained from its underlying undirected graph $G$ by replacing each edge of $G$ with the corresponding arcs of both directions, that is, $D=\overleftrightarrow{G}.$ 
A digraph $D$ is {\em minimally strong} if $D$ is strong but $D-e$ is not for every arc $e$ of $D$.


\section{Characterizations}


The following proposition can be verified using definitions of
$\kappa_k(D)$ and $\lambda_{k}(D)$.

\begin{pro}\label{proa}\cite{Sun-Yeo}
Let $D$ be a digraph of order $n$, and let $2\le k\le n$ be an integer.
Then
\begin{equation}\label{pro1}
\lambda_{k+1}(D)\leq \lambda_{k}(D) \mbox{ for every } k\le n-1
\end{equation}
\begin{equation}\label{pro2}
\kappa_k(D')\leq \kappa_k(D), \lambda_k(D')\leq \lambda_k(D) \mbox{
where $D'$ is a spanning subdigraph of $D$}
\end{equation}
\begin{equation}\label{pro3}
\kappa_k(D)\leq \lambda_k(D) \leq \min\{\delta^+(D), \delta^-(D)\}
\end{equation}
\begin{equation}\label{pro8}
\mbox{$D$ is strong if and only if $\lambda_k(D)\ge 1$.}
\end{equation}
\end{pro}

We will use the following Tillson's decomposition theorem.
\begin{thm}(Tillson's decomposition
theorem)\cite{Tillson}\label{thm02}
The arcs of $\overleftrightarrow{K}_n$ can be decomposed into
Hamiltonian cycles if and only if $n\neq 4,6$.
\end{thm}

Sun and Yeo got the following sharp bounds for $\kappa_k(D)$ and $\lambda_k(D)$.
\begin{thm}\label{thma}\cite{Sun-Yeo}
Let $D$ be a strong digraph of order $n$, and let $2\leq k\leq n$ be an
integer. Then
\begin{equation}\label{pro4}
1\leq \kappa_k(D)\leq n-1
\end{equation}
\begin{equation}\label{pro5}
1\leq \lambda_k(D)\leq n-1
\end{equation} Moreover, all bounds are sharp, and
the upper bounds hold if and only if $D\cong
\overleftrightarrow{K}_n$.
\end{thm}


By Proposition \ref{proa}(\ref{pro8}) and the fact that $\kappa_k(D)\ge 1$ if
and only if $\lambda_k(D)\ge 1$, the following directly holds:
\begin{pro}\label{pro9}
A digraph $D$ is strong if and only if $\kappa_k(D)\ge 1$ for $2\leq
k\leq n$.
\end{pro}

In the rest of this paper, we use $\mathcal{D}_{n, 1}$~(respectively, $\mathcal{D}_{n, 2})$ to denote the set of digraphs obtained from the complete digraph $\overleftrightarrow{K}_n$ by deleting an arc set $M$ such that
$\overleftrightarrow{K}_n[M]$ is a union of vertex-disjoint cycles
which cover all (respectively, $n-1$) vertices of $\overleftrightarrow{K}_n$.

\begin{lem}\label{thm2}
If $D\in \mathcal{D}_{n, 1}\cup \mathcal{D}_{n, 2}$, 
then $D$ contains $n-2$ arc-disjoint out-branchings rooted at any vertex $r$.
\end{lem}
\begin{pf}
Let $V(D)=\{u_i\mid 1\leq i\leq n\}$. Firstly, we consider the case that $D\in \mathcal{D}_{n, 1}$ such that
$\overleftrightarrow{K}_n[M]$ is a union of vertex-disjoint cycles $C_1, C_2, \ldots, C_p$
which cover all vertices of $\overleftrightarrow{K}_n$.
Without loss of generality, let $u_1$ be the root and
belong to the cycle $C_1:= u_1, u_2, \ldots, u_s, u_1$ $(s\geq 2)$. For each $3\leq i\leq s$ (if exists), let $T_{i-2}$
be an out-branching of $D$ with arc set $\{u_1u_i, u_iu_2, u_2u_{i+1}, u_iv\mid v\in V(D)\setminus \{u_1, u_2, u_i, u_{i+1}\}\}$ if $s\geq 4$; otherwise, we have $s=3$, and then let $T_1$ be an out-branching of $D$ with arc set $\{u_1u_3, u_3v\mid v\in V(D)\setminus \{u_1, u_3\}\}$.
Now consider the vertex $u_i$ for $s+1\leq i\leq n$, without loss of generality, assume that $u_i$ belongs to the cycle $C_2:= u_i, u_{i+1}, \ldots, u_{i+t-1}, u_i$~$(t\geq 2)$, where we set $u_{i+t}=u_i$ in the cycle. Let $T_{i-2}$ be an out-branching of $D$ with arc set $\{u_1u_i, u_iu_2, u_2u_{i+1}, u_iv\mid v\in V(D)\setminus \{u_1, u_2, u_i, u_{i+1}\}\}$. Hence we get $n-2$ out-branchings rooted at $u_1$ and it can be checked these out-branchings are pairwise arc-disjoint. For example, let $D\in \mathcal{D}_{5, 1}$ such that
$\overleftrightarrow{K}_5[M]$ is a union of vertex-disjoint cycles $C_1, C_2$, where $C_1:= u_1, u_2, u_3, u_1$ and $C_2:= u_4, u_5, u_4$. Let $T_1, T_2$ and $T_3$ be out-branchings rooted at $u_1$ with arc sets $\{u_1u_3, u_3u_2, u_3u_4, u_3u_5\}$, $\{u_1u_4, u_4u_2, u_2u_5, u_4u_3\}$ and $\{u_1u_5, u_5u_2, u_2u_4, u_5u_3\}$, respectively. Observe that these out-branchings are pairwise arc-disjoint.

Secondly, let $D\in \mathcal{D}_{n, 2}$ such that
$\overleftrightarrow{K}_n[M]$ is a union of vertex-disjoint cycles $C_1, C_2, \ldots, C_p$ which cover all but at most one vertex, say $u_n$, of $\overleftrightarrow{K}_n$. It suffices to consider the case that $u_n$ is the root since the argument for the remaining case 
is similar to the above paragraph. Without loss of generality, assume that the arc $u_1u_2$ belongs to one of the above cycles, say $C_1:= u_1, u_2, \ldots, u_s, u_1$ $(s\geq 2)$. Let $T_1$ be an out-branching with arc set $\{u_nu_1, u_nu_2, u_1v\mid v\in V(D)\setminus \{u_1, u_2, u_n\}\}$. For each $3\leq i\leq s$ (if exists), let $T_{i-1}$ be the out-branching of $D$ with arc set $\{u_nu_i, u_iu_2, u_2u_{i+1}, u_iv\mid v\in V(D)\setminus \{u_2, u_i, u_{i+1}, u_n\}\}$ (note that if $i=s$, then we set $u_{i+1}=u_{s+1}=u_1$ in this branching). Now consider the vertex $u_i$ for $s+1\leq i\leq n$, without loss of generality, assume that $u_i$ belongs to the cycle $C_2:= u_i, u_{i+1}, \ldots, u_{i+t-1}, u_i$ $(t\geq 2)$. Let $T_{i-1}$ be an out-branching of $D$ with arc set $\{u_nu_i, u_iu_2, u_2u_{i+1}, u_iv\mid v\in V(D)\setminus \{u_2, u_i, u_{i+1}, u_n\}\}$. Hence we get $n-2$ out-branchings rooted at $u_n$ and it can be checked that these out-branchings are pairwise arc-disjoint.
\end{pf}

\begin{lem}\label{thm3}
Let $D\in \mathcal{D}_{n, 1}\cup \mathcal{D}_{n, 2}$. For any $S\subset V(D)$ with $|S|=n-1$, $D$ contains $n-2$ internally disjoint $(S,r)$-trees, where $r\in S$.
\end{lem}
\begin{pf}
Let $V(D)=\{u_i\mid 1\leq i\leq n\}$.\\
\noindent{\bf Case 1}: 
$D\in \mathcal{D}_{n, 1}$ such that $\overleftrightarrow{K}_n[M]$ is a union of vertex-disjoint cycles $C_1, C_2, \ldots, C_p$
which cover all vertices of $\overleftrightarrow{K}_n$. Without loss of generality, let $S=V(D)\setminus \{u_2\}$ and $C_1:= u_1, u_2, \ldots, u_s, u_1$ $(s\geq 2)$.

\noindent{\bf Subcase 1.1}: $s\geq 3$. Let $D'=D[S]-\{u_1u_3\}$. Observe that $D'\in \mathcal{D}_{n-1, 1}$. By Lemma~\ref{thm2}, $D'$ contains $n-3$ arc-disjoint out-branchings $T_i~(1\leq i\leq n-3)$ rooted at any vertex $r\in S$. If $u_1$ is the root, then let $T_{n-2}$ be a tree with arc set $\{u_1u_3, u_3u_2, u_2v\mid v\in V(D)\setminus\{u_1, u_2, u_3\}\}$; if $u_i~(i\not= 1,2)$ is the root, then let $T_{n-2}$ be a tree with arc set $\{u_3u_2, u_2v\mid v\in V(D)\setminus\{u_3, u_2\}\}$ (respectively, $\{u_iu_2, u_1u_3, u_2v\mid v\in V(D)\setminus\{u_i, u_2, u_3\}\}$) when $i=3$ (respectively, $i>3$).
It can be checked that the above $n-2$ trees are pairwise internally disjoint $(S,r)$-trees for any $r\in S$ in each case.

\noindent{\bf Subcase 1.2}: $s=2$. Without loss of generality, assume that $u_3, u_4$ belong to the cycle $C_2:= u_3, u_4, \ldots, u_t, u_3$ $(t\geq 4)$.
Let $T_1$ be a tree with arc set $\{u_1u_4, u_4u_3, u_3v\mid V(D)\setminus\{u_1, u_2, u_3, u_4\}\}$; let $T_2$ be a tree with arc set $\{u_1u_3, u_3u_2, u_2v\mid V(D)\setminus\{u_1, u_2, u_3\}\}$. For $5\leq i\leq n$, 
let $T_{i-2}$ be a tree with arc set $\{u_1u_i, u_iu_4, u_4u_{i+1}, u_iv\mid v\in V(D)\setminus \{u_1, u_2, u_3, u_4, u_i\} \}$. Note that here $u_iu_{i+1}$ belongs to one of the cycles $C_2, C_3, \ldots, C_p$. It can be checked that the above $n-2$ trees are pairwise internally disjoint $(S,r)$-trees for any $r\in S$.

\noindent{\bf Case 2}: 
$D\in \mathcal{D}_{n, 2}$ such that $\overleftrightarrow{K}_n[M]$ is a union of vertex-disjoint cycles $C_1, C_2, \ldots, C_p$
which cover all but at most one vertex, say $u_n$, of $\overleftrightarrow{K}_n$. 

\noindent{\bf Subcase 2.1}: $S=V(D)\setminus \{u_n\}$. Observe that $D'=D[S]\in \mathcal{D}_{n-1, 1}$, by Lemma~\ref{thm2}, $D'$ contains $n-3$ pairwise arc-disjoint out-branchings $T_i~(1\leq i\leq n-3)$ rooted at any vertex $r\in S$. Let $T_{n-2}$ be a tree with arc set $\{ru_n, u_nv\mid v\in V(D)\setminus \{r, u_n\}\}$ where $r\in S$. It can be checked that the above $n-2$ trees are pairwise internally disjoint $(S,r)$-trees for any $r\in S$.

\noindent{\bf Subcase 2.2}: $u_n\in S$. Without loss of generality, assume that $u_2 \not\in S$ (that is, $S=V(D)\setminus \{u_2\}$) and $C_1:= u_1, u_2, \ldots, u_s, u_1$ $(s\geq 2)$. We first consider the case that $s\geq 3$. Let $D'=D[S]-\{u_1u_3\}$. Observe that $D'\in \mathcal{D}_{n-1, 2}$, by Lemma~\ref{thm2}, $D'$ contains $n-3$ arc-disjoint out-branchings $T_i~(1\leq i\leq n-3)$ rooted at any vertex $r\in S$. With a similar argument to that of Subcase 1.1, we can get $n-2$ pairwise internally disjoint $(S,r)$-trees for any $r\in S$. It remains to consider the case that $s=2$. Let $D''=D[S]-\{u_1u_n, u_nu_1\}$. Observe that $D''\in \mathcal{D}_{n-1, 1}$, by Lemma~\ref{thm2}, $D''$ contains $n-3$ pairwise arc-disjoint out-branchings $T_i~(1\leq i\leq n-3)$ rooted at any vertex $r\in S$. If $u_n$ is the root, then let $T_{n-2}$ be a tree with arc set $\{u_nu_1, u_nu_2, u_2v\mid v\in V(D)\setminus \{u_1, u_2, u_n\}\}$; otherwise, let $T_{n-2}$ be a tree with arc set $\{u_iu_n, u_nv\mid v\in V(D)\setminus \{u_1, u_i, u_n\}\}$, where $u_i~(i\not= 2, n)$ is the root. It can be checked that the above $n-2$ trees are pairwise internally disjoint $(S,r)$-trees for any $r\in S$.
\end{pf}

\begin{lem}\label{thm4}
Let $D\in \mathcal{D}_{n, 1}\cup \mathcal{D}_{n, 2}$. For any $S\subset V(D)$ with $|S|=2$, $D$ contains $n-2$ internally disjoint $(S,r)$-trees, where $r\in S$.
\end{lem}
\begin{pf}
Let $S=\{r,v\}\subset V(D)$ and $r$ be the root. If $rv\not\in M$, then let $T_1$ be the arc $rv$ and $T_u$ be an out-tree with arc set $\{ru, uv\mid uv\not\in M\}$. Otherwise, let $T_u$ be an out-tree with arc set $\{ru, uv\mid u\in V(D)\setminus \{r, v\}\}$. Observe that in both cases we get $n-2$ internally disjoint $(S,r)$-trees, as desired.
\end{pf}


We will now characterize minimally generalized $(k, \ell)$-vertex (respectively, arc)-strongly connected digraphs for some pairs of $k$ and $\ell$. 

\begin{thm}\label{thmi}
The following assertions hold:
\begin{description}
\item[(a)]~For any integer $2\leq k\leq n$, a digraph $D$ is minimally
generalized $(k, 1)$-vertex (respectively, arc)-strongly connected if and only if $D$ is minimally strong.
\item[(b)]~For any integer $2\leq k\leq n$, a digraph $D$ is minimally
generalized $(k, n-1)$-vertex (respectively, arc)-strongly connected if and only if $D\cong \overleftrightarrow{K}_n$.
\item[(c)]~For any integer $2\leq k\leq n$, a digraph $D$ is minimally
generalized $(k, n-2)$-arc-strongly connected if and only if $D\in \mathcal{D}_{n, 1}\cup \mathcal{D}_{n, 2}$; moreover, $ex'(n,k,n-2)=\mathcal{D}_{n, 1}$ and $Ex'(n,k,n-2)=\mathcal{D}_{n,
2}$.
\item[(d)]~For $k\in \{2, n-1, n\}$, a digraph $D$ is minimally generalized $(k,
n-2)$-vertex-strongly connected if and only if $D\in \mathcal{D}_{n, 1}\cup \mathcal{D}_{n, 2}$; moreover, $ex(n, k, n-2)=\mathcal{D}_{n, 1}$ and $Ex(n, k, n-2)=\mathcal{D}_{n, 2}$.
\end{description}
\end{thm}
\begin{pf} By Theorem~\ref{thma}, Propositions~\ref{proa}(\ref{pro8}) and~\ref{pro9}, and the well-known fact that every strong digraph has an out-branching rooted at any vertex, we have $(a)$ and $(b)$.

In the following argument,we just prove $(c)$ since the argument for $(d)$ is similar and simpler (by Lemmas \ref{thm2}, \ref{thm3} and \ref{thm4}).
If $D\in \mathcal{D}_{n, 1}\cup \mathcal{D}_{n, 2}$, by Proposition~\ref{proa} and Lemma~\ref{thm2}, we have $\lambda_k(D)\geq \lambda_n(D)\geq n-2$ for $2\leq k\leq n$.
For any $e\in A(D)$, $\min\{\delta^+(D-e), \delta^-(D-e)\}=n-3$, so
$\lambda_k(D-e)\leq n-3$ by Proposition~\ref{proa}(\ref{pro3}).
Thus, $D$ is minimally generalized $(k, n-2)$-arc-strongly
connected.

Let $D$ be minimally generalized $(k, n-2)$-arc-strongly
connected. By Theorem~\ref{thma}, we have $D\not
\cong \overleftrightarrow{K}_n$, that is, $D$ can be obtained from a
complete digraph $\overleftrightarrow{K}_n$ by deleting a nonempty
arc set $M$. To end our argument, we need the following proposition. Let
us start from a simple yet useful observation, which follows from
Proposition~\ref{proa}.

\begin{obs}\label{pro:HT}
No pair of arcs in $M$ has a common head or tail.
\end{obs}
\vspace{3mm}

Thus, $\overleftrightarrow{K}_n[M]$ must be a union of
vertex-disjoint cycles or paths, otherwise, there are two arcs of
$M$ such that they have a common head or tail, a contradiction with
Observation~\ref{pro:HT}.

\begin{pro}\label{pro:HT2}
$\overleftrightarrow{K}_n[M]$ does not contain a path of order at least two.
\end{pro}
\begin{pf}
Let $M'\supseteq M$ be a set of
arcs obtained from $M$ by adding some arcs from
$\overleftrightarrow{K}_n$ such that the digraph
$\overleftrightarrow{K}_n[M']$ contains no path of order at least
two. Note that $\overleftrightarrow{K}_n[M']$ is a supergraph of
$\overleftrightarrow{K}_n[M]$ and is a union of vertex-disjoint
cycles which cover all but at most one vertex of
$\overleftrightarrow{K}_n$. By Proposition~\ref{proa} and Lemma \ref{thm2}, we have $\lambda_k(\overleftrightarrow{K}_n[M'])\geq n-2$, so
$\overleftrightarrow{K}_n[M]$ is not minimally generalized $(k, n-2)$-arc-strongly
connected, a contradiction.
\end{pf}
\vspace{3mm}



It follows from Proposition~\ref{pro:HT2} and its proof that $\overleftrightarrow{K}_n[M]$ must be a union of vertex-disjoint cycles which cover all but at most one vertex of $\overleftrightarrow{K}_n$, which completes the proof.
\end{pf}

\section{The functions $f(n,k,\ell)$, $g(n,k,\ell)$, $F(n,k,\ell)$ and $G(n,k,\ell)$}

By definition, we can get the following proposition.
\begin{pro}\label{pro6} The following assertions hold:
\begin{description}
\item[(a)]~A digraph $D$ is minimally generalized $(k, \ell)$-vertex-strongly
connected if and only if $\kappa_k(D)= \ell$ and $\kappa_k(D-e)= \ell-1$ for any arc $e\in A(D)$.
\item[(b)]~A digraph $D$ is minimally generalized $(k, \ell)$-arc-strongly
connected if and only if $\lambda_k(D)= \ell$ and $\lambda_k(D-e)= \ell-1$ for any arc $e\in A(D)$.
\end{description}
\end{pro}
\begin{pf} \noindent{\bf Part $(a)$}. The direction ``if" is clear by definition, it suffices to prove the direction ``only if". Let $D$ be a minimally generalized $(k, \ell)$-vertex-strongly connected digraph. By definition, we have $\kappa_k(D)\geq \ell$ and $\kappa_k(D-e)\leq \ell-1$ for any arc $e\in A(D)$. Then for any set $S \subseteq V(D)$ with $|S|=k$, there is a set $\mathcal{D}$ of $\ell$ internally disjoint $(S, r)$-trees, where $r\in S$ is a root. As $e$ must belong to one and only one element of $\mathcal{D}$, we are done. The argument for \noindent{\bf Part $(b)$} is similar.
\end{pf}

For $2\leq k\leq n$ and $1\leq \ell \leq n-1$, let
$s(n, k, \ell)$ ($t (n, k, \ell)$, respectively) be the minimum size of a
strong digraph $D$ with order $n$ and $\kappa_k(D)=\ell$
($\lambda_k(D)=\ell$, respectively).

\begin{lem}\label{pro11}
For any $2\leq k\leq n$ and $1\leq \ell \leq n-1$,
$$s(n, k, \ell)=f(n, k, \ell), t(n, k, \ell)=g(n, k, \ell).$$
\end{lem}
\begin{pf}
Let $\mathfrak{D}(n,k,\ell)$ be the set of all strong digraphs $D$ with order $n$ and $\kappa_k(D)=\ell$. Let $$\mathfrak{D}'(n,k,\ell)=\{D\mid D\in \mathfrak{D}(n,k,\ell), |A(D)|=s(n,k,\ell)\}.$$
Recall that $\mathfrak{F}(n,k,\ell)$ is the set of all minimally generalized $(k,
\ell)$-vertex-strongly connected digraphs with order $n$, and $$ex(n,k,\ell)=\{D\mid D\in \mathfrak{F}(n,k,\ell), |A(D)|=f(n,k,\ell)\}$$ where $f(n,k,\ell)=\min\{|A(D)| \mid D\in \mathfrak{F}(n,k,\ell)\}.$

By Proposition~\ref{pro6}, $\mathfrak{F}(n,k,\ell)$ is the set of all strong digraphs $D$ with order $n$ such that $\kappa_k(D)= \ell$ and $\kappa_k(D-e)= \ell-1$ for any arc $e\in A(D)$. Hence, $\mathfrak{F}(n,k,\ell) \subseteq \mathfrak{D}(n,k,\ell)$ and so $s(n, k, \ell)\leq f(n,k,\ell)$.

Let $D\in \mathfrak{D}'(n,k,\ell)$. Then $\kappa_k(D)=\ell$ and $\kappa_k(D-e)\leq  \ell-1$ for any arc $e\in A(D)$, that is, $D\in \mathfrak{F}(n,k,\ell)$. This means that $\mathfrak{D}'(n,k,\ell)\subseteq \mathfrak{F}(n,k,\ell)$ and so $s(n, k, \ell)\geq f(n,k,\ell)$. Hence, $s(n, k, \ell)=f(n, k, \ell)$. The equality $t(n, k, \ell)=g(n, k, \ell)$ can be proved similarly.
\end{pf}

We still need the following result, see, e.g., Corollary 5.3.6 of \cite{Bang-Jensen-Gutin}.

\begin{thm}\label{2n-2-thm}
Every strong digraph $D$ on $n$ vertices has a strong spanning subgraph $H$ with at most $2n-2$ arcs and equality holds only if $H$ is a symmetric digraph whose underlying undirected graph is a tree.
\end{thm}

We will now prove our second main result:

\begin{thm}\label{thmh} The following assertions hold:
\begin{description}
\item[(a)]\begin{itemize}
\item $f(n,k,\ell)\geq n\ell$ for any two integers $2\leq k\leq n$ and $1\leq \ell \leq n-1$; moreover, the bound can be attained if $\ell=1$, or, $2\leq \ell\leq n-1$ and $k=n$.
\item $f(n,k,\ell)\leq a(\ell)+2\ell(n-\ell)$ for the case $n\geq
k+\ell$, where \[
a(\ell)=\left\{
   \begin{array}{ll}
     2{\ell\choose 2}, &\mbox {$\ell\geq 2$;}\\
     0, &\mbox {$\ell=1$.}
   \end{array}
   \right.
\]Especially, $f(n,k,\ell)\leq 2\ell(n-\ell)$ when $n\geq
k+2\ell$. \\Moreover, both bounds are sharp.
\end{itemize}
\item[(b)]~$g(n,k,\ell)= n\ell$ for any two integers $2\leq k\leq n$ and
$1\leq \ell \leq n-1$.
\item[(c)]\begin{itemize}
\item $F(n,k,\ell)\geq n\ell$ for any two integers $2\leq k\leq n$ and
$1\leq \ell \leq n-1$; especially, $F(n,k,\ell)\geq 2\ell(n-\ell)$ when $n\geq
k+2\ell$. \\Moreover, both bounds are sharp.
\item $F(n, k, \ell)=(\ell+1)(n-1)$ if $k\in \{2, n-1, n\}$ and $\ell=n-2$, or, $2\leq k\leq n$ and $\ell\in \{1, n-1\}$.
\end{itemize}
\item[(d)]\begin{itemize}
\item $G(n,k,\ell)\geq n\ell$ for any two integers $2\leq k\leq n$ and
$1\leq \ell \leq n-1$; especially, $G(n,k,\ell)\geq 2\ell(n-\ell)$ when $n\geq
k+2\ell$. \\Moreover, both bounds are sharp.
\item $G(n, k, \ell)=(\ell+1)(n-1)$ if $2\leq k\leq n$ and $\ell\in \{1, n-2, n-1\}$.
\end{itemize}
\end{description}
\end{thm}
\begin{pf}

\noindent{\bf Part $(a)$}. We first prove the lower bound.
By Proposition~\ref{proa}(\ref{pro3}), for all digraphs $D$ and $k
\geq 2$ we have $\kappa_k(D) \leq \delta^+(D)$ and $\kappa_k(D) \leq
\delta^-(D)$. Hence for each $D$ with $\kappa_k(D)=\ell$, we have
that $\delta^+(D), \delta^-(D)\geq \ell$, so $|A(D)|\geq n\ell$ and
then $f(n,k,\ell)\geq n\ell.$

We now prove the sharpness of the lower bound. For the case that $\ell=1$, let $D$ be a dicycle $\overrightarrow{C_n}$. Clearly, $D$ is minimally generalized
$(k,1)$-vertex-strongly connected, and we know $|A(D)|=n$, so
$f(n,k,1)= n$. For the case that $k=n \not\in \{4,6\}$ and $2\leq \ell\leq n-1$,
let $D\cong \overleftrightarrow{K_n}$. By Theorem \ref{thm02}, $D$
can be decomposed into $n-1$ Hamiltonian cycles $H_i~(1\leq i\leq
n-1)$. Let $D_{\ell}$ be the spanning subdigraph of $D$ with arc
set $A(D_{\ell})=\bigcup_{1\leq i\leq \ell}{A(H_i)}$. Clearly, we
have $\kappa_n(D_{\ell})\geq \ell$ for $2\leq \ell\leq n-1$.
Furthermore, by Proposition~\ref{proa}(\ref{pro3}), we have
$\kappa_n(D_{\ell})\leq \ell$ since the in-degree and out-degree of
each vertex in $D_{\ell}$ are both $\ell$. Hence,
$\kappa_n(D_{\ell})= \ell$ for $2\leq \ell\leq n-1$. For any $e\in
A(D_{\ell})$, we have
$\delta^+(D_{\ell}-e)=\delta^-(D_{\ell}-e)=\ell-1$, so
$\kappa_n(D_{\ell}-e)\leq \ell-1$ by
Proposition~\ref{proa}(\ref{pro3}).
Thus, $D_{\ell}$ is minimally generalized $(n, \ell)$-vertex-strongly
connected. As $|A(D_{\ell})|=n\ell$, we have $f(n,n,\ell)\leq
n\ell$. From the lower bound that $f(n,k,\ell)\geq n\ell$, we have
$f(n,n,\ell)= n\ell$ for the case that $2\leq \ell\leq n-1, n\not\in
\{4,6\}$. For the case that $n=6$, let $C_1$ be the cycle $u_1, u_2,
u_3, u_4, u_5, u_6, u_1$, $C_2={C_1}^{\rm rev}$, $C_3$ be the cycle
$u_1, u_3, u_5, u_2, u_6, u_4, u_1$, $C_4={C_3}^{\rm rev}$. For
$1\leq \ell\leq 4$, let $D_{\ell}$ be the digraph with vertex set
$\{u_i\mid 1\leq i\leq 6\}$ and arc set $\bigcup_{1\leq i\leq
\ell}{A(C_i)}$, let $D_{5}=\overleftrightarrow{K}_6$. It can be
checked that $D_{\ell}$ is minimally generalized $(n,
\ell)$-vertex-strongly connected. As $|A(D_{\ell})|=n\ell$, we have
$f(n,n,\ell)\leq n\ell$. From the lower bound that $f(n,k,\ell)\geq
n\ell$, we have $f(n,n,\ell)= n\ell$ for the case that $2\leq
\ell\leq n-1, n=6$. For the case that $n=4$, let $C_1$ be the cycle
$u_1, u_2, u_3, u_4, u_1$, $C_2={C_1}^{\rm rev}$. For $1\leq
\ell\leq 2$, let $D_{\ell}$ be the digraph with vertex set
$\{u_i\mid 1\leq i\leq 4\}$ and arc set $\bigcup_{1\leq i\leq
\ell}{A(C_i)}$, let $D_{3}=\overleftrightarrow{K}_4$. With a similar
but simpler argument, we can deduce that $f(n,n,\ell)= n\ell$ for
the case that $2\leq \ell\leq n-1, n=4$. Hence, $f(n,k,\ell)= n\ell$
when $2\leq \ell\leq n-1$ and $k=n$.

To prove the upper bound, we need to construct the following two digraphs $H_1$ and $H_2$ as follows:

Let $H_1$ be a symmetric digraph whose underlying undirected graph is $K_{\ell}\bigvee \overline{K}_{n-\ell}$~($n\geq k+\ell$), i.e. the graph obtained from disjoint graphs $K_{\ell}$ and $\overline{K}_{n-\ell}$ by adding all edges between the vertices in $K_{\ell}$ and $\overline{K}_{n-\ell}$. Let $V(H_1)=W_1\cup U_1$ such that $W_1=V(K_{\ell})=\{w_i\mid 1\leq i\leq \ell\}$ and
$U_1=V(\overline{K}_{n-\ell})=\{u_j\mid 1\leq j\leq n-\ell\}$.

Let $H_2=\overleftrightarrow{K}_{\ell, n-\ell}$~($n\geq k+2\ell$),
the complete bipartite digraphs with two parts $W_2$ and $U_2$, where $W_2=\{w_i\mid 1\leq i\leq \ell\}$ and $U_2=\{u_j\mid 1\leq j\leq n-\ell\}$.

\begin{pro}\label{proH1H2}~The following assertions hold:
\begin{description}
\item[(i)]~For $n\geq k+\ell$, $\kappa_k(H_1)=\ell$.
\item[(ii)]~For $n\geq k+2\ell$, $H_2$ is minimally generalized $(k,
\ell)$-vertex(arc)-strongly connected.
\end{description}
\end{pro}
\begin{pf}
For $(i)$, let $S_1$ be any $k$-subset of vertices of $V(H_1)$ such that $|S_1\cap W_1|=s$ ($s\leq \ell$) and $|S_1\cap U_1|=k-s$ ($k-s\leq n-\ell$ since $n\geq k+\ell$). Without loss of generality, let $w_i\in S_1$
for $1\leq i\leq s$ and $u_j\in S_1$ for $1\leq j\leq k-s$. For $1\leq
i\leq s$, let $T'_i$ be a tree with edge set
$$\{w_iu_1, w_iu_2, \dots , w_iu_{k-s}, u_{k-s+i}w_1, u_{k-s+i}w_2, \dots ,
u_{k-s+i}w_{s}\}.$$ For $s+1\leq j\leq \ell$, let
$T'_j$ be a tree with edge set $$\{w_ju_1, w_ju_2, \dots , w_ju_{k-s},
w_jw_1, w_jw_2, \dots , w_jw_{s}\}.$$ This is reasonable since $(k-s)+s=k\leq n-\ell$. It is not hard to obtain an $(S_1, r)$-tree $D'_i$ from $T'_i$ by adding appropriate directions to edges of $T'_i$ for any $r\in S_1$. Observe that $\{D'_i\mid 1\leq i\leq s\}\cup \{D'_j\mid s+1\leq j\leq \ell\}$ is a set of $\ell$ internally disjoint $(S_1, r)$-trees, so $\kappa_{S_1, r}(H_1)\geq \ell$, and then $\kappa_k(H_1)\geq \ell$. Combining this with the bound that
$\kappa_k(H_1)\leq \min\{\delta^+(H_1), \delta^-(H_1)\}=\ell$, we have $\kappa_k(H_1)=\ell$.

For $(ii)$, Let $S_2$ be any $k$-subset of vertices of $V(H_2)$ such that $|S_2\cap W_2|=s$ ($s\leq \ell$) and $|S_2\cap U_2|=k-s$ ($k-s< n-\ell$ since $n\geq k+2\ell$). Without loss of generality, let $w_i\in S_2$ for $1\leq i\leq s$ and $u_j\in S_2$ for $1\leq j\leq k-s$. For $1\leq i\leq s$, let $T''_i$ be a tree with edge set $$\{w_iu_1, w_iu_2, \dots , w_iu_{k-s}, u_{k-s+i}w_1, u_{k-s+i}w_2, \dots, u_{k-s+i}w_{s}\}.$$ For $s+1\leq j\leq \ell$, let $T''_j$ be a tree with edge set $$\{w_ju_1, w_ju_2, \dots , w_ju_{k-s}, w_ju_{k+j-s}, w_ju_1, w_ju_2, \dots , w_ju_{k-s}\}.$$ This is reasonable since $(k-s)+s+(\ell-s)=k+\ell-s\leq n-\ell$. It is not hard to obtain an $(S_2, r)$-tree $D''_i$ from $T''_i$ by adding appropriate directions to edges of $T''_i$ for any $r\in S_2$. Observe that $\{D''_i\mid 1\leq i\leq s\}\cup \{D''_j\mid s+1\leq j\leq \ell\}$ is a set of $\ell$ internally disjoint $(S_2, r)$-trees, so $\kappa_{S_2, r}(H_2)\geq \ell$, and then $\kappa_k(H_2)\geq \ell$. Combining this with the bound that
$\kappa_k(H_2)\leq \min\{\delta^+(H_2), \delta^-(H_2)\}=\ell$, we have $\kappa_k(H_2)=\ell$. Observe that $\kappa_k(H_2-e)\leq \min\{\delta^+(H_2-e), \delta^-(H_2-e)\}=\ell -1$. Hence, $H_2$ is minimally generalized $(k, \ell)$-vertex-strongly connected. Since the above $\ell$ internally disjoint trees are also arc-disjoint, it can be similarly proved that $H_2$ is minimally generalized $(k, \ell)$-arc-strongly connected.
\end{pf}
\vspace{2mm}

Clearly, we have $|A(H_1)|=a(\ell)+2\ell(n-\ell)$ and $|A(H_2)|=2\ell(n-\ell)$. So $f(n, k, \ell)= s(n, k, \ell)\leq a(\ell)+2\ell(n-\ell)$ when $n\geq k+\ell$ by Lemma~\ref{pro11} and Proposition~\ref{proH1H2}($i$); especially, $f(n, k, \ell)\leq 2\ell(n-\ell)$ when $n\geq k+2\ell$ by Proposition~\ref{proH1H2}($ii$).
By Theorem~\ref{thmi}($b$), when $\ell=n-1$, we have $f(n, k, \ell)= n(n-1)= a(\ell)+2\ell(n-\ell)$. Recall that when $\ell=1$, we have $f(n, k, \ell)= 2(n-1)$. Therefore, the above two bounds are sharp.

\noindent{\bf Part $(b)$}. By Proposition~\ref{proa}(\ref{pro3}),
for all digraphs $D$ and $k \geq 2$ we have $\lambda_k(D) \leq
\delta^+(D)$ and $\lambda_k(D) \leq \delta^-(D)$. Hence for each $D$
with $\lambda_k(D)=\ell$, we have that $\delta^+(D), \delta^-(D)\geq
\ell$, so $|A(D)|\geq n\ell$ and then $g(n,k,\ell)\geq n\ell.$ Now
consider the graph $D_{\ell}$ in \noindent{\bf Part $(a)$}. By
Proposition \ref{proa}(\ref{pro1}), we have $\lambda_k(D_{\ell})\geq
\lambda_n(D_{\ell})= \kappa_n(D_{\ell}) \geq \ell$ for $2\leq k\leq
n, 1\leq \ell\leq n-1$. Furthermore, by
Proposition~\ref{proa}(\ref{pro3}), we have $\lambda_k(D_{\ell})\leq
\ell$ since the in-degree and out-degree of each vertex in
$D_{\ell}$ are both $\ell$. Hence, $\lambda_k(D_{\ell})= \ell$ for
$2\leq k\leq n, 1\leq \ell\leq n-1$. For any $e\in A(D_{\ell})$, we
have $\delta^+(D_{\ell}-e)=\delta^-(D_{\ell}-e)=\ell-1$, so
$\lambda_k(D_{\ell}-e)\leq \ell-1$ by
Proposition~\ref{proa}(\ref{pro3}).
Thus, $D_{\ell}$ is minimally generalized $(k, \ell)$-arc-strongly
connected. As $|A(D_{\ell})|=n\ell$, we have $g(n,k,\ell)\leq
n\ell$. From the lower bound that $g(n,k,\ell)\geq n\ell$, we have
$g(n,k,\ell)= n\ell$.

\noindent{\bf Parts $(c)$ and $(d)$}.
By Theorem~\ref{thmi}, the following assertions hold: $F(n, k, \ell)=(\ell+1)(n-1)$ if
$k\in \{2, n-1, n\}$ and $\ell=n-2$, or, $2\leq k\leq n$ and $\ell=n-1$; $G(n, k, \ell)=(\ell+1)(n-1)$ if $2\leq k\leq n$ and $\ell\in \{n-2, n-1\}$.
Let $D$ be a minimally generalized $(k, 1)$-vertex (respectively, arc)-strongly connected digraph. By Theorems~\ref{thmi} and~\ref{2n-2-thm}, we have $|A(D)|\le 2(n-1)$ and the bound can be attained when $D$ is a symmetric digraph whose underlying undirected graph is a tree. Furthermore, we have $F(n, k ,1)=G(n, k, 1)=2(n-1)$.

By the assertions $(a)$ and $(b)$, and the fact that $f(n,k,\ell)\leq F(n,k,\ell)$, $g(n,k,\ell)\leq G(n,k,\ell)$, we directly have $F(n,k,\ell)\geq n\ell$ and $G(n,k,\ell)\geq n\ell$. Moreover, both lower bounds can be attained when $\ell=n-1$. Furthermore, we have $F(n, k, \ell)\geq 2\ell(n-\ell)$ and $G(n, k, \ell)\geq 2\ell(n-\ell)$ when $n\geq k+2\ell$ by Proposition~\ref{proH1H2}($ii$), and these two bounds can be attained for the case that $\ell=1$.
\end{pf}


\section{Discussions}

For $k\in \{2, n-1, n\}$, the minimally generalized $(k,
n-2)$-vertex-strongly connected digraphs are characterized in Theorem~\ref{thmi}. It is natural to extend this result to the case of a general $k$, 
like that of  the minimally generalized $(k, n-2)$-arc-strongly connected digraphs in Theorem~\ref{thmi}.

\begin{op}\label{op1}
Characterize the minimally generalized $(k, n-2)$-vertex-strongly connected digraphs for $2\leq k\leq n$.
\end{op}

Recall that in the proof of Theorem~\ref{thmi}, we use the monotone property of $\lambda_k$, that is, $\lambda_{k+1}(D)\leq \lambda_{k}(D)$ for every $2\leq k\le n-1$ (Proposition~\ref{proa}(\ref{pro1})). However, this property does not hold for the parameter $\kappa_k$ as shown in\cite{Sun-Yeo}, so we need to find other approach to solve Problem~\ref{op1}.

In Theorem~\ref{thmh}, we give sharp lower bounds for
$F(n,k,\ell)$ and $G(n,k,\ell)$, but we still cannot give nice upper bounds for these two functions. So It would be interesting to study the following question.

\begin{op}\label{op2}
Find sharp upper bounds for $F(n,k,\ell)$ and $G(n,k,\ell)$ for all $k\ge 2$ and
$\ell\ge 2$.
\end{op}

Let $|D_{\ell}^+|$ (respectively, $|D_{\ell}^-|)$ denote the
number of vertices of out-degree $\ell$ (in-degree
$\ell$) of a digraph $D$. By Propositions~\ref{proa}(3) and \ref{pro6}, we have $\delta^+(D), \delta^-(D)\geq \ell$ for a minimally generalized $(k, \ell)$-vertex-strongly connected digraph $D$.
It would be interesting to bound $|D_{\ell}^+|$ (respectively, $|D_{\ell}^-|)$.
Note that similar questions for minimally (strongly) connected (di)graphs were
discussed in the literature, see e.g. \cite{Halin, Kameda, Mader,
Mader2}.

\begin{op}\label{op9}
Does $|D_{\ell}^+|$ (respectively, $|D_{\ell}^-|$) $>0$ hold for every
minimally generalized $(k, \ell)$-vertex-strongly connected digraph $D$?
\end{op}

Here is a stronger question.

\begin{op}\label{op10}
Does $|D_{\ell}^+|$ (respectively, $|D_{\ell}^-|$) $\geq \ell +1$ hold
for every minimally generalized $(k, \ell)$-vertex-strongly connected digraph $D$?
\end{op}

Similar to problems \ref{op9} and \ref{op10}, the
following problems are also of interest.
\begin{op}\label{op11}
Does $|D_{\ell}^+|$ (respectively, $|D_{\ell}^-|$) $>0$ hold for every
minimally generalized $(k, \ell)$-arc-strongly connected digraph $D$?
\end{op}

\begin{op}\label{op12}
Does $|D_{\ell}^+|$ (respectively, $|D_{\ell}^-|$) $\geq \ell +1$ hold
for every minimally generalized $(k, \ell)$-arc-strongly connected digraph $D$?
\end{op}

A digraph $D$ is called {\em minimally $\ell$-vertex-strongly connected} if $D$ is $\ell$-vertex-strongly connected, but $D-e$ is not for any arc $e\in A(D)$. By definition, a 1-vertex-strongly connected digraph is also a minimally strongly connected digraph. Mader obtained the following result on minimally $\ell$-vertex-strongly connected digraphs:

\begin{thm}\label{thm03}\cite{Mader2}
For every minimally $\ell$-vertex-strongly connected digraph $D$, $|D_{\ell}^+|$ (respectively, $|D_{\ell}^-|$) $\geq \ell +1$ holds.
\end{thm}

By Theorems~\ref{thmi} and \ref{thm03}, we have the following supports for the above four problems:
Problems \ref{op11} and \ref{op12} are true for any pair of $k$ and $\ell$ with $2\leq k\leq n$ and $\ell\in \{1, n-2, n-1\}$; Problems \ref{op9} and \ref{op10} are true for any pair of $k$ and $\ell$ satisfying: $2\leq k\leq n$ and $\ell\in \{1, n-1\}$, or, $k\in \{2, n-1, n\}$ and $\ell=n-2$.

\vskip 1cm


\end{document}